\documentclass[11pt]{article}

\usepackage{amsmath,amsthm}
\usepackage{fullpage}
\usepackage{authblk}
\usepackage{lscape}
\usepackage{hyperref}

\newtheorem{prop}{Proposition}[section]
\newtheorem{remark}{Remark}[section]
\newtheorem{theorem}{Theorem}[section]

\usepackage{amssymb}
\usepackage{color}

\def\R{\mathbb{R}}
\def\N{\mathbb{N}}
\def\Z{\mathbb{Z}}
\def\E{\mathcal{E}}
\def\F{\mathcal{A}}
\def\FG{\mathcal{F}}

\def\A{\mathcal{A}}
\def\LL{\mathcal{L}}

\newcommand{\Nd}{N_{\rm c}}

\numberwithin{equation}{section}

\newcommand{\h}{{\delta}}
\newcommand{\hh}{{\delta_f}}

\title{An overview of {\sl a posteriori} error estimation and post-processing methods for nonlinear eigenvalue problems} 

\author[1]{Genevi\`eve Dusson\footnote{genevieve.dusson@math.cnrs.fr}}
\affil[1]{Universit\'e de Franche-Comt\'e,
Laboratoire de math\'ematiques de Besançon,
UMR~CNRS~6623,
16, route de Gray,
25000 Besançon,
France}
  \author[2]{Yvon Maday\footnote{yvon.maday@sorbonne-universite.fr}}
\affil[2]{Sorbonne Université, CNRS, Université Paris Cité, Laboratoire Jacques-Louis Lions (LJLL), F-75005 Paris, France et Institut Universitaire de France.}

\date{}

\begin{document}

\maketitle

\renewcommand{\thefootnote}{\fnsymbol{footnote}}

\renewcommand{\thefootnote}{\arabic{footnote}}

\begin{abstract} 
In this article, we present an overview of  different {\sl a posteriori} error analysis and post-processing methods proposed in the context of nonlinear eigenvalue problems, e.g. arising in electronic structure calculations for the calculation of the ground state and  compare them. We provide two equivalent error reconstructions based either on a second-order Taylor expansion of the minimized energy, or a first-order expansion of the nonlinear eigenvalue equation. We then show how several {\sl a posteriori} error estimations as well as post-processing methods
can be formulated as specific applications of the derived reconstructed errors, and we compare their range of applicability as well as numerical cost and precision.

\end{abstract}

\baselineskip 14pt

\setlength{\parindent}{1.5em}

\setcounter{section}{0}

\section{Introduction} \label{section1}

Nonlinear eigenvalue problems occur in many mathematical models used in science and engineering such as the calculation of the vibration modes of a mechanical structure
in the framework of nonlinear elasticity, the ground state of the Gross--Pitaevskii equation  describing the steady states of Bose--Einstein
condensates~\cite{Pitaevskii2003-ew}, or of the  Hartree--Fock and Kohn--Sham equations (\cite{lin2019numerical} used to calculate ground state electronic structures of molecular systems in quantum chemistry and materials science~(see~\cite{Cances2003-og} for a mathematical
introduction).

The approximation of the solutions to such problems is of major importance and different methods to compute them are proposed depending on the various applications at stake. 
These approximation methods are based on different ingredients. The first one is related to  the notion of degrees of freedom, associated with the basis sets used to approximate the solutions of these problems, which leads to discrete problems that can be solved, eventually, on a computer. 
Second, the resulting discrete problems are, per force, nonlinear; hence efficient algorithms must be designed to solve these problems accurately with a reasonable computational cost.
Regarding the problems studied in this article, two main families of methods exist. First, there exist algorithms directly minimizing the energy functional, such as saddle point problems solved with Newton-type methods~\cite{Bao2003-ec,Caliari2009-zk}, projected Sobolev gradient flow methods~\cite{kazemi2010minimizing,Garcia-Ripoll2001-th,Heid2021-vb, Henning2020-ha,Raza2009-mr,Zhang2022-ue,Altmann2021-ko}. The second type of methods solve the corresponding first-order Euler--Lagrange equations, in the form of a nonlinear eigenvalue problem, based on the Self-Consistent Field (SCF) algorithm~\cite{McWeeny1956-cv,Roothaan1951-ig,Cances2000-pt,Cances2000-jo,Dion2007-yt,Upadhyaya2021-wt}.
In~\cite{Cances2021-jx}, an analysis comparing direct minimization methods and SCF algorithms is proposed, see also the references therein.

These two approximation ingredients must be well tuned so that the approximate solutions are close enough to their corresponding exact ones.
Moreover the exact solutions being unknown, arguments to estimate the error between exact and approximate solutions need to be proposed. 
Such arguments often rely on so-called {\sl a priori} and {\sl a posteriori} analysis, where the {\sl a priori} analysis aims at providing convergence rates of the approximate solutions towards the exact ones, while the {\sl a posteriori} analysis derives error bounds that should only depend on quantities computable from the approximate solutions the accuracy of which one wants to qualify.
These last years, a large number of articles about the numerical analysis of approximation of solutions to problems set in the form of a partial differential equation in fields like fluid mechanics and solid mechanics have been published, see e.g. the monograph~\cite{Verfurth2013-sj} on this subject, focusing mainly on characterizing the number of degrees of freedom necessary to reach a given accuracy.

In the context of the current article, the {\sl a priori}  analysis for the class of nonlinear eigenvalue problems is  quite recent and relies on the papers \cite{Cances2003-og,Zhou2007-ci,Cances2010-zh,Chen2010-pz,Langwallner2010-fd,Cances2012-bo,Chen2011-hh} and the references therein. These articles consider the discretization error, that is the error due to the use of a given number of degrees of freedom in the approximation of the problem of interest. They provide convergence rates i.e. upper bounds for the decay rate of the error -- both on the eigenvectors and the eigenvalues -- when the number of degrees of freedom increases.

These {\sl a priori}
approaches allow first to state that: i) it is possible to achieve a satisfying approximation, provided that computing resources are sufficient; this is the notion of convergence, then, ii) for a problem in dimension $d$, multiplying the number of degrees of freedom by $2^d$ leads to an error decay by a factor $2^r$ where $r$ is related to the order of convergence of the method. Such results are nevertheless insufficient in general since most of the time, the affordable discretizations are limited due to the high computational cost of the methods 
which may be polynomial in the number of degrees of freedom e.g. $[2^d]^p$, where $p$ can be $\simeq 5$ or $10$.
Hence the interest rather goes in estimating, once a computation has been done with a given discretization, the magnitude of the error. 
This is precisely the aim for which {\sl a posteriori}  approaches (estimators and indicators) have been designed. 
These analyses do not usually give rates of convergence but instead provide actual figures that are upper bounds (and also, most of the time now, lower bounds)  of the errors between exact and approximate quantities of interest such as in this context: energies, eigenvalues, wavefunctions (eigenfunctions). 
Also they may allow for an optimization of the choice of the degrees of freedom (such as basis functions) in order to minimize the computational cost to reach a given accuracy.

As far as we know, the first article in the direction of {\sl a posteriori} estimates for nonlinear eigenvalue problems is
\cite{Maday2003-eb}, where the {\sl a posteriori} analysis of the Hartree--Fock problem was  performed and upper and lower error bounds were proposed for the ground state energy. Other contributions have been proposed since, e.g. \cite{Dai2008-fo,Chen2010-pz,Dai2015-jq}, where error bounds for nonlinear eigenvalue problems are presented in the context of finite element discretizations, see also {\sl a posteriori} results for Fourier (planewave) discretizations in \cite{Cances2014-lb,Dusson2016-ys}.

Let us also mention two-grid methods which have been designed not directly to provide error bounds, but to obtain accurate solutions at a low computational cost. 
In these methods, a first step consists in computing  a coarse solution by performing a full calculation -- in this case the resolution of a nonlinear eigenvalue problem -- with a  limited number of degrees of freedom, which should be reasonably cheap because of the small number of involved degrees of freedom. 
The second step is to perform a simpler computation, e.g. to solve a linear eigenvalue problem or a boundary value problem, with a larger number of degrees of freedom, to improve the coarse solution. 
For nonlinear eigenvalue problems, two-grid methods have been proposed for example in~\cite{Henning2014-fh,Cances2018-ow} for a Gross--Pitaevskii type equation and~\cite{Cances2016-vy} for Kohn--Sham models, the latter relying on a perturbation method. 
We will see in the following that the analysis techniques used in these two-grid methods can be very similar to the ones used in the former {\sl a posteriori} methods.

\medskip

The different references quoted above have some methodology in common but are applied to different problems, possibly with different boundary conditions, and discretized with different methods. This makes a straightforward comparison of these approaches quite difficult. In this article, we intend to compare them in a unified framework, clarifying their main similarities and differences (see Table \ref{table:tablepostprocess}). To do so, we will present all the methodologies in the framework of a simple nonlinear eigenvalue problem.
More precisely, we  place ourselves in the periodic setting
where thus
the domain $\Omega$ is the unit cell of a periodic lattice ${\cal R}$ of $\R^d$ and $
X = H^1_\#(\Omega)$,
with $d=1$, $2$ or $3$. Then for all $s \in \R^+$ and $k \in \N$\footnote{Note that $L^2_\# (\Omega)$ coincides with $L^2(\Omega)$},
\begin{eqnarray*}
H^s_\#(\Omega) & = & \left\{
v_{|\Omega}, \; v \in H^s_{\rm loc}(\R^d) \; | \; v \; \text{is} \; {\cal
  R}\mbox{-periodic} \right\}, \\
  H^{-s}_\#(\Omega) & = & [H^s_\#(\Omega)]' \quad (\hbox{dual space of} \ H^s_\#(\Omega)),  \\
C^k_\#(\Omega) & = & \left\{
v_{|\Omega}, \; v \in C^k(\R^d) \; | \; v \; \text{is}  \; {\cal
  R}\mbox{-periodic} \right\} .
\end{eqnarray*}
We then focus on a particular class of (linear and) nonlinear eigenvalue
problems arising in the study of variational models of the form
\begin{equation} \label{eq:min_pb_u}
I = \inf \left\{ \E(v), \; v \in X, \;  \int_{\Omega} v^2 = 1 \right\},
\end{equation}
associated to an energy functional $\E$  of the form\begin{equation} \label{eq:energyy}
\E(v) = \frac 12 a(v,v) + \frac \mu2 \int_\Omega G(v^2(x)) \, dx, \quad v\in X,
\end{equation}
where $\mu=1$ (nonlinear case) is the case of interest but we may also have $\mu=0$ that corresponds to the linear case, and where 
\begin{enumerate} 
\item the term in the integral involving a given function $G$ depends the square of solution ($v^2$)
 to mimic the electronic structure calculation problems we are mostly interested in such as the Hartree--Fock and Kohn--Sham problems. In these problems, the nonlinearity depends on the electronic density $\rho$, which corresponds to $v^2$ when the energy functional depends only on one function in $X$. In what follows, we denote by $g$ the derivative of $G$.

    \item the bilinear form $a$ is defined by
$$
a(u,v) = \int_\Omega (A \nabla u) \cdot \nabla v +
\int_\Omega V uv.
$$

  \item we make the following assumptions on the elements appearing in the energy
\begin{eqnarray}
\!\!\!\!\!\!\!\!\!\!\!\!\!\!\!\!\!\!\!\!\!\!\!\!\!\!\!\!\!\!\!\!\!\!\!
&\bullet& A \in (L^\infty({\Omega}))^{d\times d};\,
A(x) \mbox{ is symmetric for almost all } x \in \Omega \nonumber;\\
\!\!\!\!\!\!\!\!\!\!\!\!\!\!\!\!\!\!\!\!\!\!\!\!\!\!\!\!\!\!\!\!\!\!\!
&& \exists \,\alpha > 0 \mbox{ such that } \xi^T A(x) \xi
\ge \alpha |\xi|^2 \mbox{, } \forall \, \xi \in \mathbb{R}^d \mbox{ and almost all } x \in \Omega;
% \label{eq:Hyp1} \\
\nonumber
\\
\!\!\!\!\!\!\!\!\!\!\!\!\!\!\!\!\!\!\!\!\!\!\!\!\!\!\!\!\!\!\!\!\!\!\!
%&& \nonumber \\
%\!\!\!\!\!\!\!\!\!\!\!\!\!\!\!\!\!\!\!\!\!\!\!\!\!\!\!\!\!\!\!\!\!\!\!
&\bullet& V \in L^p(\Omega) \mbox{ for some } p > \max(1,d/2) ;
\nonumber
\\
% \label{eq:Hyp3} \\
\!\!\!\!\!\!\!\!\!\!\!\!\!\!\!\!\!\!\!\!\!\!\!\!\!\!\!\!\!\!\!\!\!\!\!
%&& \nonumber \\
%\!\!\!\!\!\!\!\!\!\!\!\!\!\!\!\!\!\!\!\!\!\!\!\!\!\!\!\!\!\!\!\!\!\!\!
&\bullet& G \in C^1([0,+\infty),\mathbb{R}) \cap C^2((0,\infty),\mathbb{R}), \,g=G', \, ~g(0)=0 \mbox{ and}~
g' > 0 \mbox{ on}~ (0,+\infty);
\label{eq:Hyp4}\mbox{} \\
\!\!\!\!\!\!\!\!\!\!\!\!\!\!\!\!\!\!\!\!\!\!\!\!\!\!\!\!\!\!\!\!\!\!\!
&& \exists \, 0 \le q < 2, \; \exists \,C \in \mathbb{R}_+ \mbox{ such that } \forall  t \ge 0, \;
|g(t)| \le C (1+t^q); 
\nonumber \\
% \label{eq:Hyp5}\\
\!\!\!\!\!\!\!\!\!\!\!\!\!\!\!\!\!\!\!\!\!\!\!\!\!\!\!\!\!\!\!\!\!\!\!
&\bullet&  g'(t)t \mbox{ is locally bounded on } [0,+\infty). 
\nonumber 
% \\
% \label{eq:Hyp6}
\end{eqnarray}
\end{enumerate}
There is no loss of generality in assuming in \eqref{eq:Hyp4} that $g(0)=0$ since the minimizers of \eqref{eq:min_pb_u} are not modified if $G(t)$ is replaced with $G(t)+ct$, due to the normalization constraint imposed on the solutions of \eqref{eq:min_pb_u}.

\medskip

\noindent
It is well known (see e.g. Lemma 2 in \cite{Cances2010-zh}) that under these assumptions, problem (\ref{eq:min_pb_u}) has exactly two minimizers $u$ and $-u$, one of them, say $u$, being positive on $\Omega$.  In all what follows, $u$ will be the positive minimizer  of \eqref{eq:min_pb_u}. Let us introduce the Fock operator, defined, for any $u\in X$ by $\F_u= D_v\E(u)$, i.e.
\begin{equation}\label{eq:Fock-def}
\F_u = - \hbox{div}A \nabla + V + \mu g(u^2).
\end{equation}
We then denote by $\A$ the operator corresponding to the linear case, i.e. for $\mu=0$, so that
\begin{equation}\label{eq:Fock-def2}
\F = - \hbox{div}A \nabla + V.
\end{equation}
Writing the Euler--Lagrange equation of problem \eqref{eq:min_pb_u}, we obtain that the function $u$ is solution to the nonlinear eigenvalue problem
\begin{equation} \label{eq:Euler-equ}
\forall v \in X, \quad \langle \F_u u-\lambda
u,v \rangle_{X',X} = 0 ,
\end{equation}
where the eigenvalue $\lambda\in \mathbb{R}$ is the Lagrange multiplier associated with the constraint $\|u\|_{L^2 }=1$. Equation~\eqref{eq:Euler-equ}, complemented with the constraint $\|u\|_{L^2 }=1$, reads
\begin{equation}
 \label{eq-VPNL}
\left\{ \begin{array}{l}
\F_u u = \lambda u, \\
\| u\|_{L^2 } = 1. \end{array} \right.
\end{equation}
Note that, for any $v\in X$, $\F_v$ is a linear self-adjoint operator on $L^2 (\Omega)$ with form domain $X$.
It can then be inferred from (\ref{eq-VPNL}) that $u\in X\cap C^0(\overline{\Omega})$, $u>0$ in $\Omega$, and $\lambda$ is the lowest eigenvalue of problem~\eqref{eq-VPNL}, called the ground state eigenvalue of $\A_u$. Note however that there exist cases where the lowest eigenvalue does not correspond to the lowest energy, such as a Gross--Pitaevskii equation with a rotating magnetic field, see~\cite[Section 6.1]{Altmann2021-ko}. Regarding quantum chemistry equations, such as the Kohn--Sham or Hartree--Fock equations, a rule called Aufbau principle states that the ground state indeed corresponds to the lowest eigenvalues of the nonlinear operator $\A_u$. An important point is that $\lambda$ is a {\em simple} eigenvalue of $\A_u$. These results are classical; their proofs are recalled in~\cite{Cances2010-zh}.

We now consider a family of finite-dimensional subspaces $(X_{\delta})_{\delta>0}$ of $X$ (conforming approximation $X_{\delta} \subset X$) such that
\begin{equation}\label{eq:assumption-space}
\forall v\in X,\quad\lim_{\delta\to 0}\min_{v_\delta\in X_\delta}\|v-v_{\delta}\|_X=0.
\end{equation}
An example is the Fourier discretization. Let assume that $\Omega = (0,2\pi)^d$. For any $k\in\Z^d$, we define the planewave $e_k$ by $e_k(x) = (2\pi)^{-d/2} e^{ik\cdot x}$. We then consider the discretization space indexed by a parameter $M$ that grows when the space grows (that is corresponding to $1/\delta$ above)
\begin{equation}\label{eq:PWSpace}
    X_M = {\rm Span} \left\{ 
    e_k, \; k\in \Z^d, \; |k| \le M
    \right\},
\end{equation}
where $|k|$ denotes the $l^2$-norm of the so-called wave-vector $k$. 

The variational approximations of \eqref{eq:min_pb_u} then consists in solving
\begin{equation}\label{dis-min-pro}
I_{\delta}=\inf\left\{\E(v_\delta), \; v_\delta\in X_\delta, \; \int_{\Omega}v_{\delta}^2=1\right\}.
\end{equation}
Problem \eqref{dis-min-pro} has at least one minimizer $u_{\delta}$ such that $(u,u_\delta)_{L^2 } \ge 0$, which satisfies
\begin{equation}
 \label{eq:discrete_nl_eig_pb}
\forall v_\delta\in X_\delta, \quad \langle \A_{u_\delta}u_\delta,v_\delta \rangle_{X',X} = \lambda_\delta (u_\delta,v_\delta)_{L^2 },
\end{equation}
for some $\lambda_\delta\in \mathbb{R}$. It is easily seen that (see e.g., \cite{Cances2010-zh, Zhou2003-vc})
\begin{equation}\label{convergence-1grid}
\lim_{\delta \to 0} \|u-u_\delta\|_X=0,
\end{equation}
or, in words, that the approximate ground state eigenfunction converges to the exact ground state eigenfunction in $H^1_\#$-norm, from which we deduce that $I_\delta$ and $\lambda_\delta$ converge to $I$ and $\lambda$, respectively, when $\delta$ goes to $0$. Optimal convergence rates have been obtained in~\cite{Cances2010-zh} (under stronger assumptions on the nonlinearity $G$) for spectral Fourier discretizations and also for finite element discretizations.
First, under suitable and realistic hypotheses on $\A$, there holds
\begin{equation}
 \label{eq:x1-2}
	\|u-u_\h\|_{X} \lesssim \min_{v_\h \in X_\h} \|u - v_\h\|_{X},
\end{equation}
where we denote by $a \lesssim b$ the inequality $a \le C b$, with $C>0$ a constant that is independent of the discretization parameter $\delta$.
Also, the eigenvalues converge much faster, similarly as in the linear case where the eigenvalues converge quadratically compared to  the eigenvectors,  i.e.
\begin{equation}
 \label{eq:x2}
	|\lambda_\h-\lambda|  \lesssim \|u-u_\h\|_{X}^2, 
 \quad \text{ for } \mu = 0.
\end{equation}
In the nonlinear case, there is an additional term involving an $L^{6/(5-2q)}$-norm
which can be absorbed in the $\|u-u_\h\|_{X}^2$-term under additional regularity assumption that we do not detail here (see e.g.~\cite[Remark 3]{Cances2010-zh} for a precise comment about this).
Moreover, the $L^2$-norm of the error $u-u_\h$ as well as its negative Sobolev norms converge faster than the $H^1_\#$-norm of the error.
Finally, the dual norm of the residual behaves like the $H^1_\#$-norm of the error, i.e.
\[
	\|u-u_\h\|_{X}\lesssim\|\A_{u_\h}u_\h-\lambda_\h u_\h\|_{X'} \lesssim \|u-u_\h\|_{X}.
\]
The first inequality is not trivial as the problem is an eigenvalue problem, but it was shown for the Gross--Pitaevskii equation in \cite{Dusson2016-ys}, and for the Laplace eigenvalue problem in \cite{Cances2017-oc}.

\medskip

The outline of this article is the following. In Section \ref{sec:rec_error}, we present two equivalent ways of approximating the error between the exact and approximate solutions, the first one relying on the minimization problem \eqref{eq:min_pb_u}, the second on the eigenvalue problem \eqref{eq:Euler-equ}. We then present different post-processing methods relying on the presented reconstruction of the error, and approximations thereof.
In Section \ref{sec:sec3}, we show that this reconstructed error also appears in {\sl a posteriori} error estimations proposed for this problem, and compare different contributions, namely exposing the main features and range of applicabilities, such as considered models and discretization methods.

\section{Two approaches for a derivation of the reconstructed error}
\label{sec:rec_error}

In this section, we present two different ways of estimating the error in this context. One is based on the energy minimization problem, and consists in looking at the second-order Taylor expansion of the energy functional. The second one relies on a first-order Taylor expansion of the nonlinear eigenvalue problem. Since the nonlinear eigenvalue problem corresponds to the first-order Euler--Lagrange equations of the minimization problem, it is natural that the two approaches are ultimately equivalent. 
However, due to the nonlinear structure of the equations and the norm constraint on the solution vector, the derivation is not straightforward, and we present the two derivations to highlight their similarities and differences, as well as pointing out the different references following these schemes.

Note that the equations obtained below for the reconstructed error are similar to what is presented in~\cite{Schmidt2020-ni} for a generic equation (without constraints). Also, in~\cite{Cances2022-uf}, such Taylor expansion and first-order error reconstruction is also proposed for the computation of ground state energies in planewave electronic structure calculations for materials systems, involving several eigenvalues.

\subsection{Approach based on the energy minimisation problem}

In this approach, presented and analyzed in~\cite{Maday2003-eb} on the Hartree--Fock problem, the initial idea is to provide a lower and an upper bound to the ground state energy $\E(u)$ from the knowledge of $\E(u_\delta)$. Of course, due to the variational statement of the conforming discretization, the following upper bound is classical
\begin{equation*}
\E(u_\delta)\ge \E(u).
\end{equation*}
To get a lower bound, the idea is to consider the second-order Taylor expansion of the energy.
Before doing so, let us introduce a notation: for any $v\in X$,
\begin{equation*}
\Lambda_v = \langle \F_v(v), v \rangle_{X',X}, 
\end{equation*}
so that the smallest eigenvalue  in \eqref{eq-VPNL}  satisfies 
$\lambda = \Lambda_{u}$. 
We then introduce the Lagrangian of the problem defined for $v\in X,$ $\nu\in \R,$ by
\begin{equation*}
\LL(v,\nu)= \E(v) - \nu \left(\int_\Omega v^2 -1\right), 
\end{equation*}
and define $\E^w(v)$ as being the value of the Lagrangian at some $\Lambda_w$, i.e.
\begin{equation*}
\E^w(v) = \LL(v,\Lambda_w)=\E(v) - \Lambda_w\left(\int_\Omega v^2 -1\right).
\end{equation*}
For $w\in X$, denoting by $D_v\E^w$ the differential of $v\in X \mapsto \E^w(v)$, there holds
\begin{equation}
\label{eq:x1}
\forall w\in X, \quad \langle [D_v\E^{u}](u), w \rangle_{X',X} = 0, 
\end{equation}
and if $u_\delta$ is a solution to the discrete problem~\eqref{eq:discrete_nl_eig_pb}, then 
\begin{equation} \label{equa2.7}
\forall w_\delta\in X_\delta, \quad \langle [D_v\E^{u_\delta}](u_\delta), w_\delta\rangle_{X',X}  = 0 .
\end{equation}
Writing the second-order Taylor expansion of the difference in energies between the two minima over $X$ and over $X_\h$ gives
 \begin{eqnarray*}
\E(u_\delta)-\E(u) & =&  \E^{u_\delta}(u_\delta)-\E^{u_\delta}(u)\\
&=& \langle [D_v\E^{u_\delta}](u_\delta), u_\delta-u\rangle_{X',X}  - \frac{1}{2}
\langle [D_v^2\E^{u_\delta}](u_\delta)(u_\delta-u),u_\delta-u\rangle_{X',X}  + o(\| u -u_\delta\|_X^2).
\end{eqnarray*}
Noting that $(u - u_\h,u_\h)_{L^2} = - \frac{1}{2} \|u_\h - u\|_{L^2}^2$, and defining $\varepsilon_0 = \|u_\h - u\|_{L^2}$ we can state
\begin{align*}
    u -u_\delta &= -  \frac{1}{2} \varepsilon_0^2 u_\h +  w, \quad \hbox{with} \quad  w\perp u_\delta, 
\end{align*}
 where $\|w\|_X \simeq \varepsilon = \| u_\delta -u \|_X$, and of course $\varepsilon_0 \le \varepsilon$. Using the above decomposition and~\eqref{equa2.7}, the energy difference can be written as
 \begin{eqnarray*}
\E(u_\delta)-\E(u)
&=& \langle [D_v \E^{u_\delta}](u_\delta),  \frac12 \varepsilon_0^2 u_\delta -  w \rangle_{X',X}  \\
&& - \frac{1}{2}
\langle [D_v^2\E^{u_\delta}](u_\delta)\left( \frac12 \varepsilon_0^2 u_\delta -  w\right) , \frac12 \varepsilon_0^2 u_\delta -  w\rangle_{X',X}  + o(\varepsilon^2)\\
&=& - \langle [D_v \E^{u_\delta}](u_\delta),    w \rangle_{X',X}  - \frac{1}{2}
\langle [D_v^2\E^{u_\delta}](u_\delta)(w),  w\rangle_{X',X}  + o(\varepsilon^2).
\end{eqnarray*}
Introducing  the bilinear form $a_v \equiv [D_v^2\E^v](v)(.,.) = \langle [D_v^2\E](v)(.),.\rangle_{X',X}  -\Lambda_v \langle.,.\rangle_{X',X} $, the previous equation reads
\begin{equation}
\label{eq:aaac}
\E(u_\delta)-\E(u)  = - \langle [D_v\E^{u_\delta}](u_\delta),    w \rangle_{X',X}  - \frac{1}{2} a_{u_\delta}(w,w) + o(\varepsilon^2).
\end{equation}
Using the fact that the lower eigenvalue is simple, that is there is a gap $\lambda_2-\lambda >0$ between the two first eigenvalues of $\F_u$,  see e.g. Lemma 1 in \cite{Cances2010-zh}, there holds:
 \begin{prop}\label{prop:z1}
There exists a constant $c_u>0$ such that, for any $v\in X$, $v\perp u$,
\begin{equation*}
a_u(v,v)\ge c_u \| v\|_{X}^2.
\end{equation*}
 \end{prop}
\noindent From this result, we obtain the following proposition.
 \begin{prop} Assume that $\varepsilon := \| u-u_\h\|_X $ is small enough, then,
there exist a constant $c_{u_\delta}>0$ such that, for any $v\in X$, $v\perp u_\delta$
\begin{equation*}
a_{u_\delta}(v,v)\ge c_{u_\delta} \| v\|_{X}^2.
\end{equation*}
 \end{prop}

Using this last proposition, we derive that there exists a unique solution $\hat w\in 
X^\perp$ (where $X^\perp$ is the orthogonal to $u_\delta$ in $X$) 
called the {\sl  reconstructed error} such that
 \begin{equation}
 \label{3K9}
a_{u_\delta}(\hat w,\psi) = -  \langle [D_v\E^{u_\delta}](u_\delta),    \psi \rangle_{X',X} \quad\forall\psi \in   X
\end{equation}
which allows to rewrite (\ref{eq:aaac}) as 
\begin{equation}\label{eq:x5ff}
\E(u)  = \E(u_\delta)-  \frac{1}{2} a_{u_\delta}(\hat w,\hat w) +  \frac{1}{2} a_{u_\delta}(w- \hat w,w- \hat w) + o(\varepsilon^2)
\end{equation}
and yields the authors in~\cite{Maday2003-eb} to the inequality
\begin{equation}\label{eq:x5}
\E(u)  \ge \E(u_\delta)-  \frac{1}{2} a_{u_\delta}(\hat w,\hat w) + o(\varepsilon^2),
\end{equation}
so that $\E(u_\delta)-  \frac{1}{2} a_{u_\delta}(\hat w,\hat w)$ is an asymptotic lower bound to the exact energy $\E(u)$.

\medskip 

To compute this reconstructed error and the associated lower bound we have to solve problem~\eqref{3K9} which, of course, cannot be done exactly and thus needs to be discretized in a larger space $X_\hh$ than $X_\delta$. This problem reads: Find $\hat w_\hh\in X_\hh^\perp$ (where $X_\hh^\perp$ is the orthogonal to $u_\delta$ in $X_\hh$) 
 such that
\begin{equation}
\label{eq:HF2gridnonlin}
	\langle(\A + 2 g'(u_\h^2) u_\h^2 + g(u_\h^2) - \lambda_\h ) \hat w_\hh ,   \psi_\hh \rangle_{X',X}  = -  \langle  (\A+ g(u_\h^2) - \lambda_\h) u_\h,    \psi_\hh \rangle_{X',X} \quad\forall\psi_\hh \in    X_\hh^\perp
	.
\end{equation}
As we shall see later this can be further refined by noticing that $w$ and $\hat w$ -- or rather its discrete representation $\hat w_\hh$ --  are $\varepsilon^2$ close.
Indeed, we first note that
\begin{equation*}
[D_v \E^{u_\delta}](u) = [D_v \E^{u_\delta}](u_\delta) + [D^2_v \E^{u_\delta}](u_\delta) (u-u_\delta) + o(\varepsilon).
\end{equation*}
Since
\begin{equation*}
[D_v \E^{u_\delta}](u) = [D_v \E^{u}](u) + (\lambda_\h-\lambda) u,
\end{equation*}
we show using~\eqref{eq:x1} and~\eqref{eq:x2} that
\begin{equation}
\label{eq:2.8}
 [D_v \E^{u_\delta}](u_\delta) + [D^2_v \E^{u_\delta}](u_\delta) (u-u_\delta) =  o(\varepsilon).
\end{equation}
Combining~\eqref{3K9} and~\eqref{eq:2.8}, we deduce  that for any $\psi\in X
$
\begin{equation}
a_{u_\delta}(u-u_\delta-\hat w,\psi) =  o(\varepsilon) \|\psi\|_{X}.
\end{equation}
Since $\| u-u_\delta - w \|_{X} = O(\varepsilon^2)$ 
, there holds
\[
    a_{u_\delta}(w-\hat w,\psi) =  o(\varepsilon) \|\psi\|_{X}, 
\]
and Proposition \ref{prop:z1} together with (\ref{prop:z1}) show that
\begin{equation}
\|w-\hat w\|_{X}   =  o(\varepsilon).
\end{equation}
Using~\eqref{eq:x5ff}, we first deduce that~\eqref{eq:x5} can be improved as an  equality
\begin{equation*}\label{eq:x5'}
\E(u)  = \E(u_\delta)-  \frac{1}{2} a_{u_\delta}(\hat w,\hat w) + o(\varepsilon^2),
\end{equation*}
and second that
$u_\delta + \hat w$ is a better approximation to $u$ than $u_\delta$. However $u_\delta + \hat w$ is yet not of norm 1, in order to cure this it remains to tune $\alpha^*$ such that $\hat u_\delta = (1-\alpha^*)u_\delta + \hat w$ is of norm 1 in $L^2(\Omega)$ (which is always possible since $\| \hat w\|_X = O(\varepsilon^2)$)\footnote{note that another normalisation can be obtained by setting $\hat u_\delta = \beta^*(u_\delta + \hat w)$}. Note that $\hat u_\delta$ is computable only from the knowledge of $u_\delta$ by inverting the linear problem~\eqref{eq:HF2gridnonlin} on a finer grid.

With this we obtain a quadratic approximation both in $X$-norm and in energy. Following~\cite{Maday2003-eb}, we write

\begin{theorem} 
\label{th:theorem}
    Let us assume that $\| u -u_\h\|_X $ is small enough, then $\hat u_\delta = (1-\alpha^*)u_\delta + \hat w$ verifies 
\begin{equation*}\label{eq:x55}
\| u-\hat u_\delta\|_{X} = o(\| u- u_\delta\|_{X}).
\end{equation*}
In addition, if $G\in C^3((0,\infty),\mathbb{R})$ then
\begin{equation*}\label{eq:x55b}
\| u-\hat u_\delta\|_{X} \lesssim \| u- u_\delta\|^2_{X},
\end{equation*}
and 
\begin{equation*}\label{eq:x56}
| \E(u)  - \E(\hat u_\delta) | \lesssim \| u- u_\delta\|^4_{X}.
\end{equation*}
\end{theorem}

\noindent The two last improved estimates follow directly from the previous analysis by changing the $o(\varepsilon^q)$ by ${\mathcal O}(\varepsilon^{q+1})$ with $q=0, 1,$ or $2$ .

\begin{remark}
    The quantities $\hat u_\h$ and $\hat w$, which can be computed with the knowledge of $u_\h$ can be  used  in practice for two different (complementary) goals. First, they can be used for a refined approximation of the solution, saying that $\hat u_\h$ is a better approximation to the exact solution $u$ than $u_\h$, and similarly for the energy using $\hat w$, that is $\E(u) \simeq \E(u_\delta)-  \frac{1}{2} a_{u_\delta}(\hat w,\hat w)$. Second, they can be used to obtain a refined error bound, as $|\E(u) - \E(u_\delta)| \simeq  \frac{1}{2} a_{u_\delta}(\hat w,\hat w)$.
    These two possible approaches will be detailed separately in Sections~\ref{sec:practical}
and~\ref{sec:sec3}.
\end{remark}

\begin{remark}
Let us mention a few works which consider the nonlinear eigenvalue problem of this type with the angle of the energy minimization, namely~\cite{kazemi2010minimizing,Faou2017-zc,Henning2020-ha,Heid2021-vb,Zhang2022-ue}. In~\cite{Heid2021-vb}, the goal is to provide an adaptive procedure for the computation of the solutions of this problem, by minimizing the energy directly, by combining gradient flow iterations and adative finite element mesh refinements. In these other works, proofs for the exponential convergence of the continuous Sobolev gradient flow are provided~\cite{kazemi2010minimizing,Henning2020-ha}, respectively for the discrete gradient flow in~\cite{Faou2017-zc,Zhang2022-ue}.
Note however that these works do not provide specific error bounds on the computed solutions.
\end{remark}

\subsection{Approach based on the nonlinear eigenvalue problem}
\label{sec:approach_nonlinear_eigpb}

In the previous section, we derived the reconstructed error equation~\eqref{eq:HF2gridnonlin} starting from the energy minimization problem~\eqref{eq:min_pb_u}. In this section, we aim at arriving at the same equation starting from the eigenvalue problem~\eqref{eq:Euler-equ}, which reads: find $(u,\lambda)$ such that $\|u\|_{L^2}=1$ and

\begin{equation} 
\label{equa:2.222}
 \FG(u,\lambda):= \F_u(u)-\lambda
u  = 0.
\end{equation}

In order to improve $(u_\delta, \lambda_\delta)$, one can think of using a single step of a Newton method in a finer discrete space $X_\hh$ as introduced in the previous section that reads: 
find $\tau_\hh\in X_\hh$, $(\tau_\hh, u_\delta)_{L^2} = 0$ and $\gamma_\hh\in \R$ such that
\begin{equation}\label{eq:2.13a}
\langle [D_v \FG](u_\delta, \lambda_\delta) (\tau_\hh),\psi_\hh \rangle_{X',X} + \langle [D_\mu \FG](u_\delta, \lambda_\delta) (\gamma_\hh),\psi_\hh \rangle_{X',X} = \langle  \FG(u_\delta, \lambda_\delta), \psi_\hh \rangle_{X',X} \quad \forall \psi_\hh\in X_\hh.
\end{equation}
Note that the solution $\tau_{\hh}$ is searched in the $L^2(\Omega)$-orthogonal complement of $u_\delta$ , instead
of $X$, 
since the problem is not well-posed on $X_\h$, hence not well-conditioned on $X$. This can also be linked to
the normalization constraint $\|u_\h\|_{L^2}^2=1.$ Indeed, this condition guarantees that the first-order equation relative to the constraint is satisfied.
Since we have
\[
	\forall v\in X, \quad \forall \mu\in \R, \quad  \FG(v,\mu) = (\A+g(v^2) -\mu)v,
\]
and the differential $D_v \FG(u_\h,\lambda_\h)$ of $\FG$ at $(u_\h,\lambda_\h)$ writes
\[
	\forall v\in X, \quad D_v \FG(u_\h, \lambda_\h)(v) = 
        (\A+2 g'(u_\h^2)u_\h^2+g(u_\h^2)-\lambda_\h) v,
\]
the single step of the Newton method in the fine grid reads in a strong form as: Find $\tau_\hh\in X_\hh$, $(\tau_\hh, u_\delta)_{L^2} = 0$ such that
\begin{equation}
\label{eq:HF2gridnonlin-bis}
	(\A + 2 g'(u_\h^2) u_\h^2 + g(u_\h^2) - \lambda_\h )
	\tau_\hh = -(\A+ g(u_\h^2) - \lambda_\h) u_\h \quad \mbox{in} \quad X_\hh^{\perp}, 
\end{equation}
which is similar to~\eqref{eq:HF2gridnonlin}.
We thus propose a norm 1 improved approximation
of $u_\h$ as 
\begin{equation*}
    \widetilde u_\h = (1-\alpha')u_\delta + \tau_\hh,
\end{equation*} 
for some $\alpha'\in\R$. Then, the difference between the post-processed energy and the approximate energy allows to estimate the error between the exact energy and the approximate one, that is
\[
    \E(u) - \E(u_\h) \simeq \E(\widetilde u_\h) - \E(u_\h).
\]
Compared to the previous {\sl a posteriori} estimation~\eqref{eq:x5}, the bounds here are not guaranteed. Nevertheless, they asymptotically match the true error, provided that the space $X_\hh$ is large enough.

\begin{remark}
     The proposed procedure in fact corresponds to a standard Newton step on the modified $\mathcal F$ functional including the normalization condition $\|u_\h\|_{L^2}^2 = 1$ as defined below in~\eqref{eq:FFF}. Hence, the quadratic convergence results of the Newton method apply in this context, so that this proposed post-processing doubles the convergence rate of the $X-$norm of the eigenfunctions and also of the energy.
\end{remark}

\subsection{Practical approximations of the reconstructed error}
\label{sec:practical}

The reconstructed error presented in~\eqref{eq:HF2gridnonlin} gives a first-order approximation of the discretization error for problem~\eqref{eq-VPNL}. 
It appears that several contributions on post-processing and  error estimation for nonlinear eigenvalue problems~\cite{Maday2003-eb,Chen2014-qu,Dusson2016-ys,Cances2014-lb,Cances2018-ow,Cances2016-vy} are based on this reconstructed error, and approximations thereof. 
We try to give an overview of these methods and compare them in this section.

The main focus of this article is on nonlinear eigenvalue problems. But naturally, post-processing methods were first developed for linear eigenvalue problems -- that are themselves actually nonlinear problems -- before nonlinear eigenvalue problems. Therefore we first describe a few methods that were developed for linear eigenvalue problems. Already in 1999 by Xu and Zhou in~\cite{Xu1999-vo}, a two-grid method has been proposed to efficiently solve eigenvalue problems. It consists in first solving an eigenvalue problem on a coarse finite element mesh, and then solving a boundary value problem on a fine mesh, in order to improve the eigenvector, then an improved eigenvalue is obtained through a Rayleigh quotient. This avoids paying the full price of solving the eigenvalue problem on the fine mesh.
Later on, in~\cite{Racheva2002-ut}, another post-processing method was presented, this time proposing only an improved eigenvalue, similarly requiring the resolution of a linear boundary value problem. 

Other works include multigrid methods, such as~\cite{Neymeyr2003-ws,Lin2014-hd,Guo2017-og},
 where the idea is to first solve an eigenvalue problem on a very coarse mesh, and then to have a family of meshes, and  improve the initial solution 
 \begin{itemize}
     \item either by solving several linear problems on finer and finer meshes~\cite{Lin2014-hd},
     \item or applying some gradient recovery operator~\cite{Guo2017-og}, 
 \end{itemize}
 correcting this way the initially found eigenvector and eigenvalue.

In the same spirit as the two-grid and multi-grid post-processing methods for linear eigenvalue problems, several methods have been proposed for nonlinear eigenvalue problems. To show how these methods relate to the reconstructed error equation~\eqref{eq:HF2gridnonlin-bis}, let us first express this equation in terms of the post-processed solution $\widetilde u_\h$ and not only the correction $\tau_\hh$.
From~\eqref{eq:HF2gridnonlin-bis} there holds
\[
   ({\A} + 2 g'(u_\h^2) u_\h^2 + g(u_\h^2) - \lambda_\h )%_{|X_\h}
	(u_\h + \tau_\hh ) = 
   ({\A} + 2 g'(u_\h^2) u_\h^2 + g(u_\h^2) - \lambda_\h )%_{|X_\h}
	u_\h
   -({\A}+ g(u_\h^2) - \lambda_\h) u_\h, 
\]
hence 
\[
   ({\A} + 2 g'(u_\h^2) u_\h^2 + g(u_\h^2) - \lambda_\h )%_{|X_\h}
	(u_\h + \tau_\hh ) 
   = 2 [g'(u_\h^2)u_\h^2 ] u_\h.
\]
The problem posed on the fine grid corresponding to the reconstructed error is therefore: Find $u_\hh \in X_\hh$ such that
\begin{equation}
   \label{eq:optimal_2grid}
   ({\A} + 2 g'(u_\h^2) u_\h^2 + g(u_\h^2) - \lambda_\h )%_{|X_\h}
	u_\hh
   = 2 [g'(u_\h^2)u_\h^2 ] u_\h.
\end{equation}

Related to this post-processing, a two-level discretization technique has been presented in~\cite{Cances2018-ow} where the authors propose three different schemes for the Gross--Pitaevskii equation. All start by solving the nonlinear eigenvalue problem in a small basis.  
Then three alternatives are proposed. The first one (scheme 1) is to solve a linear eigenvalue problem on the large basis set, fixing the nonlinearity with the coarse solution. The second one (scheme 2a) consists in solving the following boundary value problem
   \[
   ({\A} + g(u_\h^2)) u_\hh = \lambda_\h u_\h
   \quad \mbox{in} \quad X_\hh,
\]
where the term $ g(u_\h^2)$ remains on the left hand side. The third scheme (called 2b) amounts to solving a linear boundary value problem on the fine space, putting the nonlinear term on the right hand side, namely solving
\begin{equation}
   \label{eq:Henning}
   {\A} u_\hh = \lambda_\h u_\h - g(u_\h^2) u_\h
   \quad \mbox{in} \quad X_\hh.
\end{equation}
The numerical analysis of the first scheme gives the following estimates for a plane wave approximation :
\begin{align*}
   \|u - u_\hh \|_{X} &\lesssim M^{-2} \|u - u_\h \|_{X}, 
   % + \| u - u_h \|_{H^1} 
   \\
   |\lambda - \lambda_\hh| 
   + \|u - u_\hh\|_{L^2} 
   &\lesssim
   \|u - u_\h \|_{X}^2 , \\
   % + H^2 \| u - u_h \|_{H^1},
   |E(u)- E(u_\hh) | & \lesssim M^{-4} |E(u)- E(u_\h) | .
\end{align*} 

%@@@@@@@@@@@@@@@@@@@@

Related to this post-processing, a two-level discretization technique was also proposed earlier in~\cite{Henning2014-fh} in the finite element context in the form of~\eqref{eq:Henning} based on the use of quasi-orthogonality properties of a Clement type operator, hence difficult to transpose to the planewave method.
They were able to improve not only the approximation of the eigenfunctions but also the eigenvalues. The same approach based on the Localized Orthogonal Decomposition (LOD)  has been recently extended to treat both the time-dependent Gross--Pitaevskii equation and the nonlinear Gross--Pitaevskii eigenvalue problem~\cite{doding2022efficient}.

Another work proposing a post-processing technique for nonlinear eigenvalue problems is~\cite{Cances2016-vy}, see~\cite{Cances2017-oc, Dusson2020-mo} for the proofs of the estimates. 
This method is based  on a perturbative expansion of the eigenvalues and eigenvectors in order to post-process the eigenfunctions and the energy at a very low computational cost. This method was presented in the case of the Kohn--Sham equations, which is  a nonlinear eigenvalue problem, where one needs to compute a few low-lying eigenstates of the considered nonlinear operator.
The method consists in solving first the full eigenvalue problem on a small planewave space, and then to post-process the eigenvectors and eigenvalues on a larger planewave using the derived perturbation expansion. This method particularly
exploits  the diagonal structure of the Laplace operator when expressed in  planewaves, which makes the post-processing particularly cheap to perform, but makes the method difficult to generalize to different types of discretizations. 
Translated on our one-eigenpair nonlinear problem, the linear boundary value problem solved on the large discretization space aims at computing $\tau_\hh \in X_{\delta_f}$ solution to 
\begin{equation}
   \label{eq:pertKS}
   (-\Delta -\lambda_\h) \tau_\hh 
   = - ({\A} + g(u_\h^2) - \lambda_\h) u_\h, \quad \text{ in } X_{\delta_f},
\end{equation}
in a case where the operator $(-\Delta -\lambda_\h)$ is diagonal, thus only two FFTs per eigenvalue are required to compute the residual in a fine grid. 
Compared to~\eqref{eq:HF2gridnonlin-bis}, some terms in the operator on the left hand side are removed. They are actually shown to be asymptotically small compared to the Laplace operator. 
In terms of errors, one obtains that the eigenvectors and energy is improved by a factor $M^{-2}$. So  the improvement of this perturbation method is limited to $M^{-2}$ for the eigenvectors or $M^{-4}$ for the energy,
whereas, in the two-grid case, the convergence rate of the eigenfunctions can be doubled. Nevertheless it is  much less expensive.

Beyond two-grid methods, there also exist multi-grid methods which use more than two grids for computing an approximation of the solution on a fine basis. For nonlinear eigenvalue problems, several of them have been proposed, such as~\cite{Jia2016-li,Hu2018-ck,Xu2020-ro}. In these three contributions, the idea, similarly as in the two-grid case, is to first solve a nonlinear eigenvalue problem on a small discretization space, and then to post-process the solution. In the multi-grid case, this post-processing consists of several steps, often including the resolution of boundary value problems on spaces of larger and larger size, or amounting to use a multigrid technique on the larger grids, in order to even avoid solving boundary value problems on those grids. In that case, the resolution of the problems cannot be directly linked to equation~\eqref{eq:optimal_2grid}, except at the first level.

\section{Using the reconstructed error for {\sl a posteriori} error estimation}
\label{sec:sec3}

The previous section, by proposing, as in~\cite{Maday2003-eb}, a better approximation $\widetilde u_\delta$ of the exact solution $u$ using the reconstructed error equation, allows to estimate the error $u-u_\delta$ by the difference between the approximate solution $u_\delta$ and the post-processed solution, i.e.
\begin{equation*}
    u-u_\delta \simeq \widetilde u_\delta - u_\delta.
\end{equation*}
This approach is in principle valid for any type of conforming discretization.

We will see below that several works use such a post-processing step to estimate the error in the context of eigenvalue problems, sometimes in a modified way in order to obtain guaranteed error bounds. In the generic context of nonlinear problems, the contribution~\cite{Schmidt2020-ni} presents an elegant way to combine a post-processing step with guaranteed error bounds. Namely, they consider a so-called ``nonsplit residual"
which corresponds to the reconstructed error $[D {\mathcal F}(u, \lambda)]^{-1}(\mathcal F(u, \lambda))$ introduced in \eqref{eq:2.13a} where 
$\mathcal F: X\times \R \rightarrow X\times\R$ is defined by
\begin{equation}
    \label{eq:FFF}
        \mathcal F(u,\lambda) := 
    \left\{
    \begin{array}{ll}
         & \F_u(u)-\lambda u  \\
        & \displaystyle \int_\Omega u^2 - 1,
    \end{array}
\right.
\end{equation}
(with a slight change in notation for  $\mathcal F$ with respect to \eqref{equa:2.222})
and $[D {\mathcal F}(u, \lambda)]: X\times\R \rightarrow X'\times \R$ denotes the differential of $\mathcal F$ at point $(u, \lambda)$.
Then (see~\cite[Theorem 1]{Schmidt2020-ni}) considering the stability constant 
\[
    \gamma(u, \lambda) = \| [D\mathcal F (u, \lambda) ]^{-1} \|_{\mathcal L(X'\times \R,X\times \R)}, 
\]
and a local nonlinearity indicator
\[
    L(\alpha) = \sup_{(v,\mu)\in \bar{B_\alpha}(u, \lambda)} \|D\mathcal{F} (v,\mu) -
    D\mathcal{F} (u, \lambda) \|_{\mathcal L(X'\times \R,X\times \R)}, 
\]
if the validity criterion 
\[
    \alpha(u, \lambda) := 2 \gamma(u, \lambda) L(2 \varepsilon(u, \lambda) ) \le 1
\]
is satisfied, with $\varepsilon(u, \lambda) = \| [D\mathcal F (u, \lambda) ]^{-1}(\mathcal F(u, \lambda))\|$
then  problem $\mathcal F(v,\mu) = 0$ has a unique solution $(u^*, \lambda^*) \in \bar B_{2\varepsilon(u, \lambda)}(u, \lambda)$ and the error is bounded by 
\[
    \|(u^*, \lambda^*)-(u, \lambda)\|_Y \le 2 \varepsilon(u, \lambda).
\]
Therefore, this gives a generic way to estimate the error, at the price of estimating the size of $\varepsilon(u, \lambda)$ and not 
exactly solving the reconstructed error equation.

\medskip

Regarding linear eigenvalue problems, several {\sl a posteriori} error estimations have been proposed, including~\cite{Larson2000-by,Gedicke2014-dm,Li2015-wm,Horger2017-jc,Nakao2019-rg,Herbst2020-yv,Giani2021-op,Cances2022-uf,Liu2022-dk}  and references therein. Often, the estimations are based on a post-processing step, which allows to obtain a more accurate solution to the problem, more or less directly related to the error bound.  

Concerning nonlinear eigenvalue problems, there exist several works proposing error estimates for the problem of interest in this article. 
Error bounds were proposed in~\cite{Chen2011-qo}, an article presenting adaptive refinement techniques for finite element simulations of Gross--Pitaevskii type equations.
This work was later extended to the finite element simulations of Kohn--Sham equations in~\cite{Chen2014-qu}. In these two works,
the error is proved to be asymptotically bounded by an error indicator (up to a constant) involving the computation $L^2$-norm of the residual.
However, the error estimates are not fully guaranteed.

In~\cite{Dusson2016-ys} we provided an {\sl a posteriori} error estimation for the Gross--Pitaevskii equation,
based on a careful two-steps approximation.
A first coarse bound based on the analysis of the first-order Taylor expansion of the error (see Section~\ref{sec:approach_nonlinear_eigpb}) 
allows to characterize the asymptotic regime and to validate when the second proposed error bound, which is close to the real error, is guaranteed. Therefore conditions guaranteeing that the second bound is valid  can be checked in practice.
More precisely, the first bound is based on Newton--Kantorovith theorem (see e.g. \cite{Caloz1997-lb}) which similarly as in~\cite{Schmidt2020-ni} provides conditions ensuring that there exists an exact solution in the vicinity of the approximate solution, and that the error between the exact and approximate solution is bounded by
\[
	\|u-u_\h\|_{X} + |\lambda-\lambda_\h| \le 2 \|[D\FG(u_\h,\lambda_\h)]^{-1}\|_{(X',\R),(X,\R)} \|\FG(u_\h,\lambda_\h)\|_{(X',\R)}.
\]
The factor 2 in this estimation allows to absorb, in a guaranteed way,  the higher-order terms of the Taylor development.
An important part of the contribution consists of showing that the differential $D\FG$ is invertible at $(u_\h,\lambda_\h)$ and to bound the norm of its inverse. Indeed, the main part of $[D\FG(u_\h,\lambda_\h)]^{-1}$ is $(\Delta)^{-1}$, that is an isometry between $X'$ and $X$, and the remaining part in $[D\FG(u_\h,\lambda_\h)]^{-1}$ is of lower order in terms of differential operator.
In the second refined bound presented in this article, the term $[D\FG(u_\h,\lambda_\h)]^{-1} \FG(u_\h,\lambda_\h)$ is  estimated from $\Delta^{-1} \FG(u_\h,\lambda_\h)$ plus complementary terms that are shown to be negligible thanks to the first bound. Asymptotically, the following bound is obtained
\[
	\|u-u_\h\|_{X} \le \alpha  \|\FG(u_\h,\lambda_\h)\|_{(X',\R)},
\]
where $\alpha$ can be taken as close to 1 as we wish when the discretization parameter $\h$ is refined, and the second component of the residual $\FG(u_\h,\lambda_\h)$ is zero since the norm constraint is exactly satisfied.
The main drawback of this method is 
the high computational cost to obtain these bounds. Indeed, to check that the necessary assumptions are satisfied, a linear eigenvalue problem in the discrete space has to be solved, for which the lowest two eigenvalues have to be computed. Note that~\cite{Dusson2016-ys} considers planewave discretization, but, in opposition to the perturbation approach, it can  be generalized to other discretization methods such as finite elements. Note finally that the above estimate is further pushed to separate the two sources of error when the method in actually implemented and an iterative algorithm is ressorted to solving the resulting nonlinear discrete problem : the discretization error and the iteration error.

Another method has been proposed more recently for the Kohn--Sham problem~\cite{Cances2022-uf} for accurately estimating the error for quantities of interest. This work also makes use of the reconstructed error equation.
Moreover the proposed bounds are computable and accurate, in the sense that they are close to the real error. They involve a post-processing step that is quite cheap, and similar to the one of~\cite{Cances2016-vy}, i.e. they only need a few FFTs on the fine space.

In the paper~\cite{Xie2019-jw}, an {\sl a posteriori} estimation for a finite element discretization was provided. The resolution of an auxiliary boundary value problem is also necessary to obtain computable upper bound of the error. The results are asymptotic in the sense that they are valid for sufficiently small meshes, but without an {\sl a posteriori} guarantee that the mesh is fine enough so that the bounds are valid. Note that, in the frame of finite element methods, local estimators based on the use of Prager--Synge techniques allow to propose alternative global bounds that can be used to improve the precision by locally refining the finite element mesh \cite{Chen2011-hh, Chen2014-qu, Cances2017-oc}.

\section{Conclusion}
% \subsection{Main characteristics} \label{sec:characteristics}

In this article, we showed how the first-order Taylor expansion of the solved equation is related to post-processing methods as well as error estimation techniques, and summarized several works based on this especially for nonlinear eigenvalue problems.
To complete the set of examples presented above, we selected a few and  compare them in a unified way in~\autoref{table:tablepostprocess}. 
Namely, we compare the problems that are originally considered, the goals of the different contributions, as well as the computational cost of the procedure.

\begin{landscape}

   \begin{table}
   \centering
   \begin{tabular}{|c|c|c|c|c|c|}
   \hline
   & & {\bf Number} & {\bf Discreti-} &  {\bf  Generic} & \\
   {\bf Method} & {\bf Equation} & {\bf  of eigen-}  & {\bf zation} & {\bf type of} & {\bf Cost} \\
   & & {\bf values} & {\bf method} & {\bf estimates} & \\
   \hline
   % \hline
   % Linear two-grid& $-\Delta u  = \lambda u$   & $K$ lowest &  Finite  & $\|u-u_\hh\|_{H^1} \lesssim \|u-u_\h\|_{H^1}^2 $ & Boundary Value \\
   % method~\cite{Xu1999-vo}  &  &   & elements (FE)  & $|\lambda-\lambda_\hh| \lesssim  |\lambda-\lambda_\h|^2 $ &  Problem (BVP) \\
   \hline
   Nonlinear &  $(-\Delta + V + g(u^2)) u $  &  & Finite  & $\|u-u_\hh\|_{H^1} \lesssim \h^3 + \|u-u_\h\|_{H^1}$
     & Boundary \\
    two-grid & $ = \lambda u$ & 1 lowest & Elements  & $|\lambda - \lambda_\hh| + \|u-u_\hh\|_{L^2}$  &  Value  \\
   method~\cite{Henning2014-fh} & & & (FE) & $\lesssim \h^4 + \h^2 \|u-u_\hh\|_{H^1}$ & Problem (BVP) \\
   \hline
   Nonlinear &  $(-\Delta + V + g(u^2)) u $  &  & FE or  &
   (in FE) if $\|u-u_\h\|_{H^1} \lesssim \h^{-\sigma}$,  & BVP or linear \\
    two-grid & $ = \lambda u$ & 1 lowest & planewaves & then  $\|u-u_\hh\|_{H^1} \lesssim \h^{-\sigma-2} + \hh^{-\sigma}$ &  eigenvalue  \\
   method~\cite{Cances2018-ow} & & & (PW) & $|E(u)-E(u_\hh)| \lesssim \h^{-2\sigma-4} + \hh^{-2\sigma}$ & problem \\
   \hline
   % \hline
   % Linear &   $(-\Delta+V) u  = \lambda u$ & $K$ lowest & PW & if $\Nd$ cutoff in momentum space,  & Residual  \\
   %  perturbation  & & & & $\|u-u_\hh\|_{H^1} \lesssim \Nd^{-2} \|u-u_\h\|_{H^1}$ & computation \\
   % method~\cite{Cances2014-lb} & & & & $|\lambda-\lambda_\hh| \lesssim \Nd^{-2} |\lambda-\lambda_\h|$  & \\
    % \hline
   Nonlinear  & $(-\Delta+V+ V_{\rho_\Phi}) \phi_i  $ & $K$ lowest & PW &  if $\Nd$  cutoff in momentum space, & Residual \\
   perturbation  & $= \lambda_i \phi_i$ & & & $\|u - u_\hh\|_{H^1} \lesssim M^{-2} \|u - u_\h\|_{H^1}$  & computation \\
   method~\cite{Cances2016-vy} & with $\Phi = (\phi_1,\ldots,\phi_K)$. &  & & Energy:  &  \\
    & & &  & $|E(u)-E(u_\hh)| \lesssim M^{-2} |E(u)-E(u_\h)|$  & \\
   \hline
   \hline
   {\sl A posteriori} estimation & $(-\Delta + V + (\rho_\Phi \star \frac{1}{|x|}) \phi_i $ & $K$ lowest & Any  &  $\|u - u_\hh\|_{H^1} \lesssim  \|u - u_\h\|_{H^1}^2$ & \\
   for Hartree--Fock~\cite{Maday2003-eb} &  $= \lambda_i \phi_i$, & &  & $|E(u)-E(u_\hh)| \lesssim  |E(u)-E(u_\h)|^2$ & BVP  \\
   &  with $\Phi = (\phi_1,\ldots,\phi_K)$. & &   &  & \\
   \hline
   {\sl A posteriori} estimation & $(-\Delta + V + g(u^2)) u  $ & $1$ lowest & PW  &  $\|u - u_\h\|_{X} \lesssim \|\mathcal F(u_\h,\lambda_\h)\|_{X'}$ & \\
   for Gross--Pitaevskii~\cite{Dusson2016-ys} &  $= \lambda u$, & &  &  & BVP  \\
   &   & &   &  & \\
   \hline
   {\sl A posteriori} estimation & $F (u) = 0  $ & no  & Any  &  $\|u - u_\h\|_{X} \lesssim \|\mathcal F(u_\h,\lambda_\h)\|_{X'}$ & \\
   for nonlinear problem~\cite{Schmidt2020-ni} &  & eigenvalue &  &  & BVP  \\
   &   & &   &  & \\
   \hline
   {\sl A posteriori} estimation & $(-\Delta + V_{\rho_\Phi}) \phi_i $ & $K$ lowest & PW  &  $\|\gamma_\Phi - \gamma_{\Phi_\hh}\|_X \lesssim  \|\mathcal F(\Phi)\|_{X'} $ & \\ 
   for Kohn--Sham~\cite{Cances2022-uf} &  $= \lambda_i \phi_i$, & &  & (density matrix) & BVP  \\
   &  with $\Phi = (\phi_1,\ldots,\phi_K)$. & &   &  & \\
   \hline
   \end{tabular}
   \caption{Comparative table of different post-processing and error estimation methods.}
   \label{table:tablepostprocess}
   \end{table}
   
\end{landscape}

\section*{Dedication}
We dedicate this article to the late Professor Roland Glowinski who was an inspiration to so many generations of applied and interdisciplinary mathematicians. The second author (YM) had the chance to follow his DEA (then master) courses at the Université Pierre et Marie Curie (then Sorbonne Universit\'e) and to benefit from his inspiring teachings both in these courses and in the lectures he was given. Even if the nonlinear eigenvalue problems he was interested in \cite{glowinski2020numerical} were of a different nature, applications to 
molecular dynamics and electronic structure calculation \cite{zhang1994formation, zhang2000distance} has been of interest for him for many years and will undoubtedly benefit from his ideas for a long time to come.

\section*{Funding}
For this project
GD is supported by
%  CALSIMLAB and the ANR within the Investissements d'Avenir programme under reference  ANR-11-IDEX-0004-02. 
%Project supported by 
the French “Investissements d'Avenir” program, project ISITE-BFC (contract ANR-15-IDEX-0003) and YM is supported by the European Research Council (ERC) under the European Union’s Horizon 2020 research and innovation program (grant agreement No 810367), project EMC2.

\bibliographystyle{siam}
\bibliography{biblio.bib}

\begin{thebibliography}{10}

\bibitem{Altmann2021-ko}
{\sc R.~Altmann, P.~Henning, and D.~Peterseim}, {\em The j-method for the
  {Gross--Pitaevskii} eigenvalue problem}, Numer. Math., 148 (2021),
  pp.~575--610.

\bibitem{Bao2003-ec}
{\sc W.~Bao and W.~Tang}, {\em Ground-state solution of {Bose--Einstein}
  condensate by directly minimizing the energy functional}, J. Comput. Phys.,
  187 (2003), pp.~230--254.

\bibitem{Caliari2009-zk}
{\sc M.~Caliari, A.~Ostermann, S.~Rainer, and M.~Thalhammer}, {\em A
  minimisation approach for computing the ground state of {Gross--Pitaevskii}
  systems}, J. Comput. Phys., 228 (2009), pp.~349--360.

\bibitem{Caloz1997-lb}
{\sc G.~Caloz and J.~Rappaz}, {\em Numerical analysis for nonlinear and
  bifurcation problems}, Handb. Numer. Anal., 5 (1997), pp.~487--637.

\bibitem{Cances2000-pt}
{\sc E.~Canc{\`e}s}, {\em {SCF} algorithms for {HF} electronic calculations},
  in Mathematical Models and Methods for Ab Initio Quantum Chemistry,
  M.~Defranceschi and C.~Le~Bris, eds., Springer Berlin Heidelberg, Berlin,
  Heidelberg, 2000, pp.~17--43.

\bibitem{Cances2018-ow}
{\sc E.~Canc{\`e}s, R.~Chakir, L.~He, and Y.~Maday}, {\em Two-grid methods for
  a class of nonlinear elliptic eigenvalue problems}, IMA J. Numer. Anal., 38
  (2018), pp.~605--645.

\bibitem{Cances2010-zh}
{\sc E.~Canc{\`e}s, R.~Chakir, and Y.~Maday}, {\em Numerical analysis of
  nonlinear eigenvalue problems}, J. Sci. Comput., 45 (2010), pp.~90--117.

\bibitem{Cances2012-bo}
{\sc E.~Canc{\`e}s, R.~Chakir, and Y.~Maday}, {\em Numerical analysis of the
  planewave discretization of some orbital-free and {Kohn-Sham} models}, Esaim
  Math. Model. Numer. Anal., 46 (2012), pp.~341--388.

\bibitem{Cances2003-og}
{\sc E.~Canc{\`e}s, M.~Defranceschi, W.~Kutzelnigg, C.~Le~Bris, and Y.~Maday},
  {\em Computational quantum chemistry: A primer}, in Handbook of Numerical
  Analysis, vol.~10, Elsevier, 2003, pp.~3--270.

\bibitem{Cances2022-uf}
{\sc E.~Canc{\`e}s, G.~Dusson, G.~Kemlin, and A.~Levitt}, {\em Practical error
  bounds for properties in {Plane-Wave} electronic structure calculations},
  SIAM J. Sci. Comput., 44 (2022), pp.~B1312--B1340.

\bibitem{Cances2014-lb}
{\sc E.~Canc{\`e}s, G.~Dusson, Y.~Maday, B.~Stamm, and M.~Vohral{\'\i}k}, {\em
  A perturbation-method-based a posteriori estimator for the planewave
  discretization of nonlinear schr{\"o}dinger equations}, C. R. Math., 352
  (2014), pp.~941--946.

\bibitem{Cances2016-vy}
\leavevmode\vrule height 2pt depth -1.6pt width 23pt, {\em A
  perturbation-method-based post-processing for the planewave discretization of
  {Kohn--Sham} models}, J. Comput. Phys., 307 (2016), pp.~446--459.

\bibitem{Cances2017-oc}
\leavevmode\vrule height 2pt depth -1.6pt width 23pt, {\em Guaranteed and
  robust a posteriori bounds for laplace eigenvalues and eigenvectors:
  Conforming approximations}, SIAM J. Numer. Anal., 55 (2017), pp.~2228--2254.

\bibitem{Cances2021-jx}
{\sc {\'E}.~Canc{\`e}s, G.~Kemlin, and A.~Levitt}, {\em Convergence analysis of
  direct minimization and {Self-Consistent} iterations}, SIAM J. Matrix Anal.
  Appl., 42 (2021), pp.~243--274.

\bibitem{Cances2000-jo}
{\sc E.~Canc{\`e}s and C.~Le~Bris}, {\em Can we outperform the {DIIS} approach
  for electronic structure calculations?}, Int. J. Quantum Chem., 79 (2000),
  pp.~82--90.

\bibitem{Chen2014-qu}
{\sc H.~Chen, X.~Dai, X.~Gong, L.~He, and A.~Zhou}, {\em Adaptive finite
  element approximations for {Kohn--Sham} models}, Multiscale Model. Simul., 12
  (2014), pp.~1828--1869.

\bibitem{Chen2011-hh}
{\sc H.~Chen, X.~Gong, L.~He, and A.~Zhou}, {\em Adaptive finite element
  approximations for a class of nonlinear eigenvalue problems in quantum
  physics}, Adv. Appl. Math. Mech., 3 (2011), pp.~493--518.

\bibitem{Chen2010-pz}
{\sc H.~Chen, X.~Gong, and A.~Zhou}, {\em Numerical approximations of a
  nonlinear eigenvalue problem and applications to a density functional model},
  Math. Methods Appl. Sci., 33 (2010), pp.~1723--1742.

\bibitem{Chen2011-qo}
{\sc H.~Chen, L.~He, and A.~Zhou}, {\em Finite element approximations of
  nonlinear eigenvalue problems in quantum physics}, Comput. Methods Appl.
  Mech. Engrg., 200 (2011), pp.~1846--1865.

\bibitem{Dai2015-jq}
{\sc X.~Dai, L.~He, and A.~Zhou}, {\em Convergence and quasi-optimal complexity
  of adaptive finite element computations for multiple eigenvalues}, IMA J.
  Numer. Anal., 35 (2015), pp.~1934--1977.

\bibitem{Dai2008-fo}
{\sc X.~Dai, J.~Xu, and A.~Zhou}, {\em Convergence and optimal complexity of
  adaptive finite element eigenvalue computations}, Numer. Math., 110 (2008),
  pp.~313--355.

\bibitem{Dion2007-yt}
{\sc C.~M. Dion and E.~Canc{\`e}s}, {\em Ground state of the time-independent
  {Gross--Pitaevskii} equation}, Comput. Phys. Commun., 177 (2007),
  pp.~787--798.

\bibitem{doding2022efficient}
{\sc C.~D{\"o}ding, P.~Henning, and J.~W{\"a}rneg{\aa}rd}, {\em An efficient
  two level approach for simulating bose-einstein condensates}, arXiv preprint
  arXiv:2212.07392,  (2022).

\bibitem{Dusson2020-mo}
{\sc G.~Dusson}, {\em Post-processing of the plane-wave approximation of
  schr{\"o}dinger equations. part {II}: {Kohn--Sham} models}, IMA J. Numer.
  Anal., 41 (2020), pp.~2456--2487.

\bibitem{Dusson2016-ys}
{\sc G.~Dusson and Y.~Maday}, {\em A posteriori analysis of a nonlinear
  {Gross--Pitaevskii-type} eigenvalue problem}, IMA J. Numer. Anal.,  (2016),
  p.~drw001.

\bibitem{Faou2017-zc}
{\sc E.~Faou and T.~J{\'e}z{\'e}quel}, {\em Convergence of a normalized
  gradient algorithm for computing ground states}, IMA J. Numer. Anal., 38
  (2017), pp.~360--376.

\bibitem{Garcia-Ripoll2001-th}
{\sc J.~J. Garc{\'\i}a-Ripoll and V.~M. P{\'e}rez-Garc{\'\i}a}, {\em Optimizing
  schr{\"o}dinger functionals using sobolev gradients: Applications to quantum
  mechanics and nonlinear optics}, SIAM J. Sci. Comput., 23 (2001),
  pp.~1316--1334.

\bibitem{Gedicke2014-dm}
{\sc J.~Gedicke and C.~Carstensen}, {\em A posteriori error estimators for
  convection--diffusion eigenvalue problems}, Comput. Methods Appl. Mech. Eng.,
  268 (2014), pp.~160--177.

\bibitem{Giani2021-op}
{\sc S.~Giani, L.~Grubi{\v s}i{\'c}, H.~Hakula, and J.~S. Ovall}, {\em A
  posteriori error estimates for elliptic eigenvalue problems using auxiliary
  subspace techniques}, J. Sci. Comput., 88 (2021).

\bibitem{glowinski2020numerical}
{\sc R.~Glowinski, S.~Leung, H.~Liu, and J.~Qian}, {\em On the numerical
  solution of nonlinear eigenvalue problems for the monge-amp{\`e}re operator},
  ESAIM: Control, Optimisation and Calculus of Variations, 26 (2020), p.~118.

\bibitem{Guo2017-og}
{\sc H.~Guo, Z.~Zhang, and R.~Zhao}, {\em Superconvergent two-grid methods for
  elliptic eigenvalue problems}, J. Sci. Comput., 70 (2017), pp.~125--148.

\bibitem{Heid2021-vb}
{\sc P.~Heid, B.~Stamm, and T.~P. Wihler}, {\em Gradient flow finite element
  discretizations with energy-based adaptivity for the {Gross-Pitaevskii}
  equation}, J. Comput. Phys., 436 (2021), p.~110165.

\bibitem{Henning2014-fh}
{\sc P.~Henning, A.~M{\aa}lqvist, and D.~Peterseim}, {\em Two-level
  discretization techniques for ground state computations of {Bose-Einstein}
  condensates}, SIAM J. Numer. Anal., 52 (2014), pp.~1525--1550.

\bibitem{Henning2020-ha}
{\sc P.~Henning and D.~Peterseim}, {\em Sobolev gradient flow for the
  {Gross--Pitaevskii} eigenvalue problem: Global convergence and computational
  efficiency}, SIAM J. Numer. Anal., 58 (2020), pp.~1744--1772.

\bibitem{Herbst2020-yv}
{\sc M.~F. Herbst, A.~Levitt, and E.~Canc{\`e}s}, {\em A posteriori error
  estimation for the non-self-consistent {Kohn--Sham} equations}, Faraday
  Discuss., 224 (2020), pp.~227--246.

\bibitem{Horger2017-jc}
{\sc T.~Horger, B.~Wohlmuth, and T.~Dickopf}, {\em Simultaneous reduced basis
  approximation of parameterized elliptic eigenvalue problems}, Esaim Math.
  Model. Numer. Anal., 51 (2017), pp.~443--465.

\bibitem{Hu2018-ck}
{\sc G.~Hu, H.~Xie, and F.~Xu}, {\em A multilevel correction adaptive finite
  element method for {Kohn--Sham} equation}, J. Comput. Phys., 355 (2018),
  pp.~436--449.

\bibitem{Jia2016-li}
{\sc S.~Jia, H.~Xie, M.~Xie, and F.~Xu}, {\em A full multigrid method for
  nonlinear eigenvalue problems}, Sci. China Math., 59 (2016), pp.~2037--2048.

\bibitem{kazemi2010minimizing}
{\sc P.~Kazemi and M.~Eckart}, {\em Minimizing the gross-pitaevskii energy
  functional with the sobolev gradient—analytical and numerical results},
  International Journal of Computational Methods, 7 (2010), pp.~453--475.

\bibitem{Langwallner2010-fd}
{\sc B.~Langwallner, C.~Ortner, and E.~S{\"u}li}, {\em Existence and
  convergence results for the galerkin approximation of an electronic density
  functional}, Math. Models Methods Appl. Sci., 20 (2010), pp.~2237--2265.

\bibitem{Larson2000-by}
{\sc M.~G. Larson}, {\em A posteriori and a priori error analysis for finite
  element approximations of self-adjoint elliptic eigenvalue problems}, SIAM J.
  Numer. Anal., 38 (2000), pp.~608--625.

\bibitem{Li2015-wm}
{\sc H.~Li and J.~S. Ovall}, {\em A posteriori eigenvalue error estimation for
  the schr{\"o}dinger operator with the inverse square potential}, Mathematics
  and Statistics Faculty Publications and Presentations,  (2015).

\bibitem{lin2019numerical}
{\sc L.~Lin, J.~Lu, and L.~Ying}, {\em Numerical methods for kohn--sham density
  functional theory}, Acta Numerica, 28 (2019), pp.~405--539.

\bibitem{Lin2014-hd}
{\sc Q.~Lin and H.~Xie}, {\em A multi-level correction scheme for eigenvalue
  problems}, Math. Comput., 84 (2014), pp.~71--88.

\bibitem{Liu2022-dk}
{\sc X.~Liu and T.~Vejchodsk{\'y}}, {\em Fully computable a posteriori error
  bounds for eigenfunctions}, Numer. Math.,  (2022).

\bibitem{Maday2003-eb}
{\sc Y.~Maday and G.~Turinici}, {\em Error bars and quadratically convergent
  methods for the numerical simulation of the {Hartree-Fock} equations}, Numer.
  Math., 94 (2003), pp.~739--770.

\bibitem{McWeeny1956-cv}
{\sc R.~McWeeny}, {\em The density matrix in self-consistent field theory i.
  iterative construction of the density matrix}, Proc. R. Soc. Lond. A Math.
  Phys. Sci., 235 (1956), pp.~496--509.

\bibitem{Nakao2019-rg}
{\sc M.~T. Nakao, M.~Plum, and Y.~Watanabe}, {\em Numerical Verification
  Methods and {Computer-Assisted} Proofs for Partial Differential Equations},
  Springer, Singapore, 2019.

\bibitem{Neymeyr2003-ws}
{\sc {Neymeyr, K. }}, {\em Solving mesh eigenproblems with multigrid
  efficiency}, Methods for Scientific Computing. Variational problems and
  Applications,  (2003).

\bibitem{Pitaevskii2003-ew}
{\sc L.~P. Pitaevskii, S.~Stringari, and S.~Stringari}, {\em {Bose-Einstein}
  Condensation}, Clarendon Press, Apr. 2003.

\bibitem{Racheva2002-ut}
{\sc M.~R. Racheva and A.~B. Andreev}, {\em Superconvergence postprocessing for
  eigenvalues}, Comput. Methods Appl. Math., 2 (2002), pp.~171--185.

\bibitem{Raza2009-mr}
{\sc N.~Raza, S.~Sial, S.~S. Siddiqi, and T.~Lookman}, {\em Energy minimization
  related to the nonlinear schr{\"o}dinger equation}, J. Comput. Phys., 228
  (2009), pp.~2572--2577.

\bibitem{Roothaan1951-ig}
{\sc C.~C.~J. Roothaan}, {\em New developments in molecular orbital theory},
  Rev. Mod. Phys., 23 (1951), pp.~69--89.

\bibitem{Schmidt2020-ni}
{\sc A.~Schmidt, D.~Wittwar, and B.~Haasdonk}, {\em Rigorous and effective
  a-posteriori error bounds for nonlinear problems---application to {RB}
  methods}, Adv. Comput. Math., 46 (2020), p.~32.

\bibitem{Upadhyaya2021-wt}
{\sc P.~Upadhyaya, E.~Jarlebring, and E.~H. Rubensson}, {\em A density matrix
  approach to the convergence of the self-consistent field iteration}, Numer.
  Algebra Control Optim., 11 (2021), p.~99.

\bibitem{Verfurth2013-sj}
{\sc R.~Verfurth}, {\em A Posteriori Error Estimation Techniques for Finite
  Element Methods (Numerical Mathematics and Scientific Computation)}, Oxford
  University Press, 1~ed., May 2013.

\bibitem{Xie2019-jw}
{\sc H.~Xie and M.~Xie}, {\em Computable error estimates for ground state
  solution of {Bose--Einstein} condensates}, J. Sci. Comput., 81 (2019),
  pp.~1072--1087.

\bibitem{Xu2020-ro}
{\sc F.~Xu and Q.~Huang}, {\em Cascadic adaptive finite element method for
  nonlinear eigenvalue problem based on complementary approach}, J. Comput.
  Appl. Math., 372 (2020), p.~112720.

\bibitem{Xu1999-vo}
{\sc J.~Xu and A.~Zhou}, {\em A two-grid discretization scheme for eigenvalue
  problems}, Math. Comput., 70 (1999), pp.~17--26.

\bibitem{zhang1994formation}
{\sc F.~Zhang, F.~Spiegelmann, E.~Suraud, V.~Frayss{\'e}, R.~Poteau,
  R.~Glowinski, and F.~Chatelin}, {\em On the formation of transient (na19) 2
  and (na20) 2 cluster dimers from molecular dynamics simulations}, Physics
  Letters A, 193 (1994), pp.~75--81.

\bibitem{zhang2000distance}
{\sc F.-S. Zhang, F.~Wang, E.~Suraud, and R.~Glowinski}, {\em A distance
  dependent tight-binding molecular dynamics model to the collision and
  thermodynamical properties of nan}, Progress of Theoretical Physics
  Supplement, 138 (2000), pp.~72--77.

\bibitem{Zhang2022-ue}
{\sc Z.~Zhang}, {\em Exponential convergence of sobolev gradient descent for a
  class of nonlinear eigenproblems}, Commun. Math. Sci., 20 (2022),
  pp.~377--403.

\bibitem{Zhou2003-vc}
{\sc A.~Zhou}, {\em An analysis of finite-dimensional approximations for the
  ground state solution of {Bose--Einstein} condensates}, Nonlinearity, 17
  (2003), p.~541.

\bibitem{Zhou2007-ci}
\leavevmode\vrule height 2pt depth -1.6pt width 23pt, {\em Finite dimensional
  approximations for the electronic ground state solution of a molecular
  system}, Math. Methods Appl. Sci., 30 (2007), pp.~429--447.

\end{thebibliography}

\end{document}